\documentclass[12pt]{article}

\usepackage{amsmath,amsthm,amsfonts,amssymb,amscd}

\allowdisplaybreaks[4]
\newtheorem {Lemma}{Lemma}[section]
\newtheorem {Theorem} {Theorem}[section]
\newtheorem {Remark}{Remark}[section]

\newtheorem {Corollary}{Corollary}[section]

\newenvironment {Proof} {\noindent {\bf Proof.}}{\quad $\square$\par\vspace{3mm}}
\newtheorem {Claim}{Claim}[section]

\begin{document}

\title{Laplacian and signless Laplacian spectral radii of graphs with fixed domination number}

\author{Rundan Xing, Bo Zhou\footnote{Corresponding author. E-mail: zhoubo@scnu.edu.cn}\\
Department of Mathematics,
South China Normal University, \\
Guangzhou 510631, P. R. China}

\date{}
\maketitle

\begin{abstract}
In this paper, we determine the maximal Laplacian and signless
Laplacian spectral radii for graphs with fixed number of vertices
and domination number, and characterize the extremal graphs
respectively.
\end{abstract}

\section{Introduction}

We consider simple undirected graphs. Let $G$ be a graph with vertex
set $V(G)$ and edge set $E(G)$. For $u\in V(G)$, let $N_G(u)$ be the
set of neighbors of vertex $u$ in $G$. The degree of vertex $u$ in
$G$, denoted by $d_G(u)$, is the cardinality of $N_G(u)$.

Let $V(G)=\{v_1,\dots,v_n\}$. The degree matrix of $G$ is the
$n\times n$ diagonal matrix $D(G)$ with its $(i,i)$-entry equal to
$d_G(v_i)$. The adjacency matrix of $G$ is the $n\times n$ matrix
$A(G)=(a_{ij})$ where $a_{ij}=1$  if $v_iv_j\in E(G)$ and $0$
otherwise. Then $L(G)=D(G)-A(G)$ is the Laplacian matrix of $G$ and
$Q(G)=D(G)+A(G)$ is the signless Laplacian matrix of $G$. Obviously,
both $L(G)$ and $Q(G)$ are all symmetric. The Laplacian spectral
radius and signless Laplacian spectral radius of $G$, denoted by
$\mu(G)$ and $q(G)$, are the largest eigenvalues of $A(G)$, $L(G)$
and $Q(G)$, respectively.

A dominating set of $G$ is a vertex subset $S$ of $G$ such that
each vertex of $V(G)\setminus S$ is adjacent to at least one vertex
of $S$. The domination number of $G$, denoted by $\gamma(G)$, is the
minimal cardinality of a dominating set of $G$. A dominating set $S$
of $G$ is said to be minimal if $|S|=\gamma(G)$.

If $G$ is an $n$-vertex graph with domination number $n$, then $G$
is the $n$-vertex empty graph, of which the adjacency, Laplacian and
signless Laplacian spectral radii are all equal to zero.

For $1\le \gamma\le n-1$, let $\mathcal{G}_{n,\gamma}$ be the set of
graphs with $n$ vertices and domination number $\gamma$.
Stevanovi\'c et al.~\cite{SAH} determined the unique graphs with
maximal spectral radius for graphs in $\mathcal{G}_{n,\gamma}$.

Recall that if $G$ contains no isolated vertices, then $\gamma(G)\le
\lfloor\frac{n}{2}\rfloor$~\cite{Ore}. Brand and Seifter \cite{BS}
gave an upper bound for Laplacian spectral radius of connected
graphs in $\mathcal{G}_{n,\gamma}$, where
$1\le\gamma\le\lfloor\frac{n}{2}\rfloor$.

In this paper, we determine the maximal Laplacian and signless
Laplacian spectral radii for graphs in $\mathcal{G}_{n,\gamma}$, and
characterize the extremal graphs respectively.

\section{Preliminaries}

Let $G$ be a graph. For $E\subseteq E(G)$, let $G-E$ be the graph
obtained from $G$ by deleting all edges of $E$. Let $\overline{G}$
be the complement of $G$. For $F\subseteq E(\overline{G})$, let
$G+F$ be the graph obtained from $G$ by adding all edges of $F$. If
$E=\{e\}$ or $F=\{f\}$, then write $G-e$ or $G+f$ instead.

\begin{Lemma}\label{G+e}\cite{He,Mo,CRS}
Let $G$ be a graph, and let $e\in E(\overline{G})$. Then
$\mu(G+e)\ge\mu(G)$.
\end{Lemma}

\begin{Lemma}\label{comp_bipar}\cite{Merris}
Let $G$ be an $n$-vertex graph. Then $\mu(G)\le n$ with equality if
and only if $\overline{G}$ is disconnected.
\end{Lemma}

%

Let $G$ be a graph. Let $\Delta(G)$ and $\overline{d}(G)$ be the
maximal degree and the average degree of $G$, respectively.
Obviously, $\overline{d}(G)=\frac{2|E(G)|}{|V(G)|}$.

\begin{Lemma}\label{max-deg}\cite{CRS}
Let $G$ be a graph. Then
\[
2\overline{d}(G)\le q(G)\le 2\Delta(G)
\]
with either equality when $G$ is connected if and only if $G$ is regular.
\end{Lemma}


Let $K_n$ be the complete graph on $n$ vertices. For an $n$-vertex
connected graph $G$, from the previous lemma, $q(G)\le 2(n-1)$ with
equality if and only if $G\cong K_n$.

Let $K_{m,n}$ be the complete bipartite graph with partite sizes $m$
and $n$, where $1\le m\le n$.

Let $G$ and $H$ be two vertex-disjoint graphs. Denote by $G\cup H$
the vertex-disjoint union of $G$ and $H$. For integer $r\ge
0$, let $rG$ be the vertex-disjoint union of $r$ copies of $G$.

An independent set of graph $G$ is a vertex subset of $G$, in
which no two vertices are adjacent in $G$.

\section{Maximal Laplacian spectral radius of graphs in $\mathcal{G}_{n,\gamma}$}

For a connected graph $G$ in $\mathcal{G}_{n,\gamma}$ with $1\le
\gamma\le \lfloor\frac{n}{2}\rfloor$, Brand and Seifter~\cite{BS}
showed that if $\gamma=1$, then $\mu(G)=n$, if $\gamma=2$, then no
better bound than  $\mu(G)\le n$ exists, and if $\gamma\ge 3$, then
$\mu(G)<n-\lceil\frac{\gamma-2}{2}\rceil$.

In this section, we determine the maximal Laplacian spectral radius
of graphs in $\mathcal{G}_{n,\gamma}$, and characterize the extremal
graphs, where $1\le\gamma\le n-1$. As a corollary, we give the
corresponding result for bipartite graphs in
$\mathcal{G}_{n,\gamma}$.

\begin{Remark}
Let $G\in\mathcal{G}_{n,1}$. Then  $d_G(u)=n-1$ for
some $u\in V(G)$, i.e.,  $K_{1,n-1}$ is  a subgraph of $G$, and
thus by Lemma~\ref{comp_bipar}, we have $\mu(G)=n$.
\end{Remark}


Let $G$ be a graph with  $U\subseteq V(G)$. For $u\in V(G)$, let
$N_G(u:U)$ be the set of neighbors of $u$ in $U$. Let
$d_U(u)=|N_G(u:U)|$. Obviously, $N_G(u:V(G))=N_G(u)$.

Let $G$ be a bipartite graph with bipartition $(U,W)$. Let $G^+$ be
the set of graphs $H$ such that $V(H)=V(G)$ and $E(G)\subseteq E(H)
\subseteq E(G)\cup E_U\cup E_W$, where $E_U=\{uv:u,v\in U\mbox{ and
}N_G(u:W)=N_G(v:W)\}$ and $E_W=\{uv:u,v\in W\mbox{ and
}N_G(u:U)=N_G(v:U)\}$. Let $\mathcal{S}^+$ be the union of all
$G^+$, where $G$ is a bipartite semi-regular graph.

Yu et al. gave the following result, where the upper bound in the
lemma was first proposed by Das \cite{Das}.

\begin{Lemma}\label{upper}\cite{YLT}
Let $G$ be a connected graph. Then
\[
\mu(G)\le\max_{uv\in E(G)}|N_G(u)\cup N_G(v)|
\]
with equality if and only if $G\in \mathcal{S}^+$.
\end{Lemma}

For $n\ge 4$, let
$\mathcal{B}_n=\{K_{a,n-a}:2\le a\le \lfloor\frac{n}{2}\rfloor\}$.

For a graph $G$ with $uv\in E(G)$, let $D_{uv}=(V(G)\setminus
(N_G(u)\cup N_G(v)))\cup\{u,v\}$.

\begin{Theorem}\label{L-dom}
Let $G\in\mathcal{G}_{n,\gamma}$, where $2\le\gamma\le n-1$. Then
$\mu(G)\le n-\gamma+2$ with equality if and only if $G\cong
H\cup(\gamma-2)K_1$, where $H\in B^+$,
$B\in\mathcal{B}_{n-\gamma+2}$ and $d_G(u)\le n-\gamma$ for $u\in
V(G)$.
\end{Theorem}

\begin{Proof}
Obviously, $\mu(G)=\mu(G_1)$ for some nontrivial connected component
$G_1$ of $G$. For $uv\in E(G_1)$, it is easily seen that $D_{uv}$
is a dominating set of $G$, implying that $\gamma\le
|D_{uv}|=n-|N_G(u)\cup N_G(v)|+2$, i.e.,
\[
|N_G(u)\cup N_G(v)|\le n-\gamma+2
\]
with equality if and only if  $D_{uv}$ is a
minimal dominating set of $G$. By Lemma~\ref{upper}, we have
\begin{eqnarray*}
   \mu(G)
&=&
   \mu(G_1)\\
&\le&
   \max_{uv\in E(G_1)}|N_{G_1}(u)\cup N_{G_1}(v)|\\
&=&
   \max_{uv\in E(G_1)}|N_G(u)\cup N_G(v)|\\
&\le&
   n-\gamma+2
\end{eqnarray*}
with equalities if and only if $G_1\in \mathcal{S}^+$ and $D_{uv}$
is a minimal dominating set of $G$ for some $uv\in E(G_1)$.

If $G\cong H\cup(\gamma-2)K_1$, where $H\in B^+$ with
$B\in\mathcal{B}_{n-\gamma+2}$ and $d_G(u)\le n-\gamma$ for
$u\in V(G)$, then by Lemmas~\ref{G+e} and~\ref{comp_bipar}, we have
$n-\gamma+2=\mu(B)\le\mu(H)\le\mu(K_{n-\gamma+2})=n-\gamma+2$,
implying that $\mu(G)=\mu(H)=n-\gamma+2$.

Suppose that $\mu(G)=n-\gamma+2$. We are to show that $G\cong
H\cup(\gamma-2)K_1$, where $H\in B^+$,
$B\in\mathcal{B}_{n-\gamma+2}$ and $d_G(u)\le n-\gamma$ for $u\in
V(G)$.

Since $G_1\in \mathcal{S}^+$, there exists a bipartite semi-regular
graph $B$ such that $G_1\in B^+$. Let $(U,W)$ be the  bipartition of
$B$. Recall that  $D_{uv}$ is a minimal dominating set of $G$ for
some $uv\in E(G_1)$. If there exist $w_1,w_2\in
D_{uv}\setminus\{u,v\}$ such that $w_1w_2\in E(G)$, then
$D_{uv}\setminus\{w_1\}$ is a dominating set of $G$ with cardinality
less than $|D_{uv}|=\gamma$, a contradiction. Thus
$D_{uv}\setminus\{u,v\}$ is an independent set of $G$.

\begin{Claim}\label{isolated}
All vertices of $D_{uv}\setminus\{u,v\}$ are isolated in $G$.
\end{Claim}

Suppose that $w\in D_{uv}\setminus\{u,v\}$ is a non-isolated vertex
of $G$. Then $w\in V(G_1)$. Suppose that $u,v$ lie in different
bipartite sets of $B$, say $u\in U$ and $v\in W$. Suppose that
$W\setminus N_{G_1}(u)\ne \emptyset$, say $a\in W\setminus
N_{G_1}(u)$. If $a\in N_{G_1}(v)$, then $av\in E_W$, and thus by the
definition of graphs in $B^+$, $N_{G_1}(a:U)=N_{G_1}(v:U)$, implying
that $a\in N_{G_1}(u)$, a contradiction. Thus $a\not\in
N_{G_1}(u)\cup N_{G_1}(v)=N_G(u)\cup N_G(v)$, implying that $a\in
D_{uv}\setminus\{u,v\}$. Since $B$ is bipartite semi-regular and
$a\in W$, $d_B(a)=d_B(v)=|N_{G_1}(v:U)|$. Since
$D_{uv}\setminus\{u,v\}$ is an independent set of $G$,
$N_{G_1}(a:U)=N_{G_1}(v:U)$, implying that $au\in E(G_1)$, a
contradiction. Thus $W\setminus N_{G_1}(u)=\emptyset$, implying that
$W=N_{G_1}(u:W)$. Similarly, $U=N_{G_1}(v:U)$. Since each vertex of
$D_{uv}\setminus\{u,v\}$ is adjacent to neither $u$ nor $v$ in $G$,
$w$ dose not exist, a contradiction. Thus  $u,v$ lie in the same
bipartite set of $B$, say $u,v\in U$. Then $uv\in E_U$, and thus by
the definition of graphs in $B^+$, $N_{G_1}(u:W)=N_{G_1}(v:W)$. Let
$N=N_{G_1}(u:W)$. Obviously, $N\ne\emptyset$, and for $a\in
N_{G_1}(u)\setminus N$ or  $a\in N_{G_1}(v)\setminus N$,
$N_{G_1}(a:W)=N$. It follows that the subgraph of $B$ induced by
$N_{G_1}(u)\cup N_{G_1}(v)=N_G(u)\cup N_G(v)$ is complete bipartite
with bipartition $((N_G(u)\cup N_G(v))\setminus N,N)$. Suppose that
$W\setminus N\ne \emptyset$, say $b\in W\setminus N$. Then $b\in
D_{uv}\setminus\{u,v\}$. Since $B$ is bipartite semi-regular,
$d_B(b)\ge |(N_G(u)\cup N_G(v))\setminus N|$. Since
$D_{uv}\setminus\{u,v\}$ is an independent set of $G$, $b$ is
adjacent to each vertex of $(N_G(u)\cup N_G(v))\setminus N$,
implying that $b\in N$, a contradiction. Thus $W=N$ and $w\in
U\setminus\{u,v\}$. Since $d_B(w)=d_B(u)$, $w$ is adjacent to each
vertex of $W$. Then $(D_{uv}\setminus\{v,w\})\cup\{a\}$ for some
$a\in W$ is a dominating set of $G$ with cardinality less than
$|D_{uv}|=\gamma$, also a contradiction. Thus each vertex of
$D_{uv}\setminus\{u,v\}$ is isolated in $G$, which proves
Claim~\ref{isolated}.

Recall that $G_1\in B^+$, where $B$ is a bipartite semi-regular graph with  bipartition $(U,W)$.

\begin{Claim}\label{com-bip}
$B$ is a complete bipartite graph.
\end{Claim}

Suppose first that $u,v$ lie in different bipartite sets of $B$, say
$u\in U$ and $v\in W$. If there exists $w\in U\setminus N_{G_1}(v)$,
then $uw\in E_U$, implying that $N_{G_1}(w:W)=N_{G_1}(u:W)$, and
thus $w\in N_{G_1}(v)$, a contradiction. It follows that $U\setminus
N_{G_1}(v)=\emptyset$, implying that $U=N_{G_1}(v:U)$. Similarly,
$W=N_{G_1}(u:W)$. Since $B$ is bipartite semi-regular, $B$ is
complete bipartite. Now suppose that $u,v$ lie in the same bipartite
set of $B$, say $u,v\in U$. Then $uv\in E_U$, and thus by the
definition of graphs in $B^+$, $N_{G_1}(u:W)=N_{G_1}(v:W)$. For
$w\in U$, $w\in N_G(u)\cup N_G(v)$, implying that $uw\in E_U$ or
$vw\in E_U$, and then $N_{G_1}(w:W)=N_{G_1}(u:W)=N_{G_1}(v:W)$. Thus
$W=N_{G_1}(u:W)$, implying that $B$ is complete bipartite. This
proves Claim~\ref{com-bip}.

If $|U|=1$ ($|W|=1$, respectively), then
$(D_{uv}\setminus\{u,v\})\cup U$ ($(D_{uv}\setminus\{u,v\})\cup W$,
respectively) is a dominating set of $G$ with cardinality less than
$|D_{uv}|=\gamma$, a contradiction. Thus $|U|,|W|\ge 2$. By
Claims~\ref{isolated} and~\ref{com-bip},  $G\cong
H\cup(\gamma-2)K_1$, where $H\in B^+$ and
$B\in\mathcal{B}_{n-\gamma+2}$.

Note that $|V(H)|=n-\gamma+2$. If
there exists a vertex $w$ of degree $n-\gamma+1$ in $G$, then $w\in
V(H)$ and $(D_{uv}\setminus\{u,v\})\cup\{w\}$ is a dominating set
of $G$ with cardinality less than $\gamma$, a contradiction. Thus each
vertex of $G$ is of degree at most $n-\gamma$ in $G$.
\end{Proof}

From the previous theorem, we easily obtain the following result for
bipartite graphs.

\begin{Corollary}\label{bipartL}
Let $G$ be a bipartite graph in $\mathcal{G}_{n,\gamma}$, where
$2\le\gamma\le n-1$. Then $\mu(G)\le n-\gamma+2$ with equality if
and only if $G\cong H\cup(\gamma-2)K_1$, where $H\in
\mathcal{B}_{n-\gamma+2}$.
\end{Corollary}

\section{Maximal signless Laplacian spectral radius of graphs in $\mathcal{G}_{n,\gamma}$}

In this section, we determine the maximal signless Laplacian
spectral radius of graphs in $\mathcal{G}_{n,\gamma}$, and
characterize the extremal graphs, where $1\le\gamma\le n-1$.

Let $G$ be a graph. For $u\in V(G)$, let $D_u=V(G)\setminus N_G(u)$.

\begin{Theorem}\label{Q-dom}
Let $G\in\mathcal{G}_{n,\gamma}$, where $1\le\gamma\le n-1$. Then
$q(G)\le 2(n-\gamma)$ with equality if and only if $G\cong
K_{n-\gamma+1}\cup (\gamma-1)K_1$ or when $\gamma\ge 2$ and
$n-\gamma$ is even, $G\cong \overline{\frac{n-\gamma+2}{2}K_2}\cup
(\gamma-2)K_1$.
\end{Theorem}

\begin{Proof}
Let $G$ be a graph with maximal signless Laplacian spectral radius
among graphs in $\mathcal{G}_{n,\gamma}$. Obviously, $q(G)=q(G_1)$
for some nontrivial connected component $G_1$ of $G$. For $u\in
V(G)$ with $d_G(u)=\Delta(G)$, it is easily seen that $D_u$ is a
dominating set of $G$, and then $\gamma\le |D_u|=n-\Delta(G)$,
implying that $\Delta(G)\le n-\gamma$ with equality if and only if
$D_u$ is a minimal dominating set of $G$. By Lemma~\ref{max-deg}, we
have
\[
q(G)=q(G_1)\le 2\Delta(G_1)\le 2\Delta(G)\le 2(n-\gamma)
\]
with equalities if and only if $G_1$ is regular and
$\Delta(G_1)=\Delta(G)=n-\gamma$, i.e., $G_1$ is
$(n-\gamma)$-regular and $D_u$ is a minimal dominating set of $G$
for some $u\in V(G_1)$. If $\gamma=1$, then it is easily seen that
$G\cong K_n$. Suppose in the following that $2\le\gamma\le n-1$.

Suppose that $q(G)=2(n-\gamma)$. Then $G_1$ is $(n-\gamma)$-regular
and for some $u\in V(G_1)$, $D_u$ is a minimal dominating set of
$G$. If there exist $v_1,v_2\in D_u\setminus\{u\}$ such that
$v_1v_2\in E(G)$, then $D_u\setminus\{v_1\}$ is a dominating set of
$G$ with cardinality less than $|D_u|=\gamma$, a contradiction. Thus
$D_u\setminus\{u\}$ is an independent set of $G$, implying that each
connected component different from $G_1$ is an isolated vertex.

Suppose that $d_{G_1}(u)\le |V(G_1)|-3$. Then there exists
$\{v,w\}\subseteq V(G_1)$ such that $uv,uw\not\in E(G_1)$. Since
$G_1$ is $(n-\gamma)$-regular, $|N_{G_1}(u)|=n-\gamma$, and
$D_u\setminus\{u\}$ is an independent set of $G$, $v$ and $w$ are
both adjacent to each vertex of $N_{G_1}(u)$, implying that
$(D_u\setminus\{v,w\})\cup\{a\}$ for some $a\in N_{G_1}(u)$ is a
dominating set of $G$ with cardinality less than $|D_u|=\gamma$, a
contradiction. Thus $d_{G_1}(u)=|V(G_1)|-1$ or $|V(G_1)|-2$.

If $d_{G_1}(u)=|V(G_1)|-1$, then since $G_1$ is
$(n-\gamma)$-regular, we have $|V(G_1)|=n-\gamma+1$ and $G_1\cong
K_{n-\gamma+1}$, implying that $G\cong K_{n-\gamma+1}\cup
(\gamma-1)K_1$.

Suppose that $d_{G_1}(u)=|V(G_1)|-2$. Then
$V(G_1)=N_G(u)\cup\{u,v\}$, where $v$ is the unique vertex in
$V(G_1)\setminus\{u\}$ which is nonadjacent to $u$. Since $G_1$ is
$(n-\gamma)$-regular, $v$ is adjacent to each vertex of
$N_{G_1}(u)$, and for $w\in N_{G_1}(u)$, $w$ is nonadjacent to
exactly one vertex of $N_{G_1}(u)\setminus\{w\}$ in $G_1$, implying
that $|N_{G_1}(u)|=n-\gamma$ is even. Thus $G_1\cong
\overline{\frac{n-\gamma+2}{2}K_2}$, which implies that $G\cong
\overline{\frac{n-\gamma+2}{2}K_2}\cup (\gamma-2)K_1$.

Conversely, if $G\cong K_{n-\gamma+1}\cup (\gamma-1)K_1$ or when
$n-\gamma$ is even, $G\cong\overline{\frac{n-\gamma+2}{2}K_2}\cup
(\gamma-2)K_1$, then by Lemma~\ref{max-deg},
$q(G)=q(G_1)=2\Delta(G_1)=2(n-\gamma)$.
\end{Proof}


If $G$ is a bipartite graph, then $L(G)$ and $Q(G)$ are unitarily
similar~\cite{Grone}. For a bipartite graph
$G\in\mathcal{G}_{n,\gamma}$ with $2\le\gamma\le n-1$, by
Corollary~\ref{bipartL}, $q(G)\le n-\gamma+2$ with equality if and
only if $G\cong H\cup(\gamma-2)K_1$, where $H\in
\mathcal{B}_{n-\gamma+2}$.

\vspace{4mm}

\noindent {\bf Acknowledgement.} \noindent {\bf Acknowledgement.}
This work was supported by the National Natural Science Foundation
of China (No.~11071089) and the Specialized Research Fund for the
Doctoral Program of Higher Education of China (No.~20124407110002).

\end{document}